\documentclass{proc-l}

\usepackage{amssymb,amsmath,amsthm,hyperref}

\newtheorem{theorem}{Theorem}[section]
\newtheorem*{main}{Main Theorem}
\newtheorem{lemma}[theorem]{Lemma}
\newtheorem{corollary}[theorem]{Corollary}

\theoremstyle{definition}
\newtheorem{definition}[theorem]{Definition}

\theoremstyle{remark}
\newtheorem{remark}[theorem]{Remark}

\DeclareMathOperator{\artanh}{artanh}
\DeclareMathOperator{\arsinh}{arsinh}

\numberwithin{equation}{section}

\begin{document}

\title[Practical {D}ehn parental bounds]
      {Practical bounds for a {D}ehn parental test}


%    author one information
\author[Robert C. Haraway III] {Robert C. Haraway III}
\address{Department of Mathematics\\
  Oklahoma State University\\
  Stillwater, OK 74078}
\email{robert.haraway@okstate.edu}
\thanks{This research supported in part by NSF grant~DMS-1006553.}

%    author two information

%    \subjclass is required.
\subjclass[2010]{Primary: 57M50; Secondary: 57-04}

\date{\today}

%    "Communicated by" -- provide editor's name; required.
\commby{David Futer}

%    Abstract is required.
\begin{abstract}
Hodgson and Kerckhoff proved a powerful theorem, half of which
they used to make Thurston's Dehn surgery theorem effective.
The calculations derived here use both halves of Hodgson and
Kerckhoff's theorem to give bounds leading towards a practical
algorithm to tell, given two orientable complete
hyperbolic 3-manifolds $M,N$ of finite volume,
whether or not $N$ is a Dehn filling of $M.$
\end{abstract}

\maketitle


\section{Introduction}
Given any two orientable complete hyperbolic 3-manifolds
$M,N$ of finite volume,
how can one determine whether or not $N$ is a Dehn filling of $M$?
To show that $N$ is a Dehn filling of $M$ is straightforward;
one can just give a Dehn filling coefficient
$c$ on $\partial M$, and show that $M$ filled along $c$ is
isometric to $N$. But it is unclear \emph{a priori} how to
show $N$ is not a Dehn filling of $M$, apart from
basic sanity checks like number of boundary components.
The following theorem, and the main result of this paper,
reduces this problem to computations of basic properties 
of hyperbolic 3-manifolds---volume,
cusp geometry, systole length, and isometry---and
avoids the need to drill geodesics out.

\begin{main}
Let $M,N$ be orientable 3-manifolds
admitting complete hyperbolic metrics
of finite volume on their interiors.
Let $\Delta V = Vol(M) - Vol(N)$.

$N$ is a Dehn filling of $M$ if and
only if either
\begin{itemize}
\item $N$ is a Dehn filling of $M$
along a coefficient of normalized length
at most $8.5$, or
\item $N$ both \emph{has} a closed simple geodesic
of length less than $0.97 \cdot \Delta V$, and
\emph{is} a Dehn filling of $M$ along
a coefficient of normalized length $\mathbf{L}$
such that
\[
\frac{8.95}{\Delta V} < \mathbf{L}^2 < \frac{13.4}{\Delta V}.
\]
\end{itemize}
\end{main}

We begin with a survey of relevant results on
the geometry of Dehn filling in section \ref{sec:DF}.
We prove the main theorem in section \ref{sec:thm};
the proof is an elementary elaboration of the
work of Hodgson and Kerckhoff in \cite{HK08}.
Finally, we discuss some issues that
might arise in using these bounds to implement
a Dehn parental test, and we conclude
with connections to prior work and ongoing projects.

\section{Hyperbolic Dehn filling}\label{sec:DF}
\emph{Hyperbolic space}, or $\mathbf{H}^3$,
is a connected, simply-connected
Riemannian 3-manifold of constant sectional curvature -1.
Such a metric space is unique up to isometry.
A \emph{hyperbolic 3-manifold}
is, for the purposes of this paper, an orientable 3-manifold
admitting a complete Riemannian metric
of finite volume that is locally isometric to $\mathbf{H}^3$.
Such a metric is, surprisingly, unique up to isometry.

Any such manifold $M$ is the interior
of a unique compact 3-manifold $\overline{M}$,
every boundary component of which is a torus.
Given such a boundary component $k$,
a \emph{slope} (on $k$) is an equivalence class
of essential simple closed curve on $k$,
two such curves being equivalent
if they are isotopic in $k$.
One can parametrize such classes as elements of
$P H_1(k;\mathbb{Q}) \approx P^1 \mathbb{Q}$,
hence the name ``slope.''

A \emph{Dehn filling coefficient} of $M$ is,
for each boundary component $k$ of $\overline{M}$,
a single choice $c(k)$ of slope on $k$,
or possibly no choice at all.
Given such a $c$, consider the following construction.
Let $K_c$ be the set of components $k$ of $\partial \overline{M}$
such that $c(k)$ is a slope on $k$.
For each $k \in K_c$, represent $c(k)$
by an essential simple closed curve $\gamma_k$.
Pick homeomorphisms $\phi_k: k \to \partial (D^2 \times S^1)$
such that $\phi_k(\gamma_k) = (\partial D^2, 1).$
Let $T =  \bigsqcup_{k \in K_c} (D^2 \times S^1)$ and
let $\Phi: \partial \overline{M} \to \partial T$
be the disjoint union of the $\phi_k$.
Then the interior of $\overline{M} \sqcup_\Phi T$
is a compact oriented 3-manifold $M'$.
The homeomorphism class $M_c$ of $M'$
is independent of our choices of $\gamma$ and $\phi_k$.
We say $M_c$ is the \emph{Dehn filling} of $M$ \emph{along $c$}.

The following theorem of Thurston inaugurated
the study of Dehn filling geometry:
\begin{theorem}[\cite{Th82}, 2.6]
  Let $L \subset N$ be a link in a 3-manifold $N$ such that
  $M = N \smallsetminus L$ is a hyperbolic 3-manifold.

  There is a finite set $S$ of slopes
  on components of $\partial \overline{M}$
  such that all Dehn fillings
  along coefficients excluding slopes from $S$
  are hyperbolic 3-manifolds.
\end{theorem}

One way to quantify this result is to provide
some measure of length for Dehn filling slopes,
then to say that all slopes of large length
(for a suitable definition of ``large'')
yield hyperbolic Dehn fillings.
The most natural choice of length is defined as follows.

Suppose $M$ is a cusped hyperbolic 3-manifold of finite volume,
and suppose $s$ is a simple closed curve
on a component $K$ of $\partial M$.
There is a maximal embedded open horoball neighborhood $B$ of $K$.
Take the completion of $B$
under the shortest-path metric on $B$.
The boundary of this completion
is a Euclidean surface homeomorphic to $K$,
which we identify with $K$.
Define \emph{the length of $s$} to be
the length of any (Euclidean) geodesic in $K$ isotopic to $s$.

The earliest quantification of Thurston's theorem
using this notion of length was the $2\pi$-Theorem,
due to Gromov and Thurston himself.
The state of the art along these lines
is the following theorem, due independently
to Agol and to Lackenby\footnote{And, in truth, also to Perelman;
Agol and Lackenby prove a condition called word-hyperbolicity, which
is equivalent to hyperbolicity by Perelman's
proof of geometrization.}:
\begin{theorem}[6-Theorem, \cite{Ag00}, \cite{La00}]
  Every Dehn filling along a slope of length at least 6 is hyperbolic.
\end{theorem}

Having proven a Dehn filling to be hyperbolic, 
one naturally may ask to what extent
its geometry differs from a ``parent''---that is,
a manifold of which it is a Dehn filling.
Thurston also showed that high-order Dehn fillings
are Gromov-Hausdorff close to the parent manifold.
We may quantify this closeness
by calculating how close geometric properties of these manifolds are.

The most important geometric property
of a hyperbolic 3-manifold is its volume.
So we are interested in the difference in volume 
$\Delta V(c) = vol(M) - vol(M_c)$
between a hyperbolic manifold $M$
and a Dehn filling $M_c$ of $M$.
Several results in the literature relate $\Delta V$
to the lengths of the Dehn filling slopes,
notably the asymptotics of Neumann and Zagier \cite{NZ85},
and the bound of Agol and Dunfield \cite{AST},
which yielded the best lower bound
on the volume of a closed orientable hyperbolic 3-manifold
until the minimum such volume
was determined outright in \cite{GMMJAMS},
\cite{GMM11}, and \cite{Milley}.
However, the bounds in \cite{NZ85} and \cite{AST}
do not serve our purposes.
We now move to bounds lying closer to what we need.

First, the following bound of Futer, Kalfagianni, and Purcell
is a necessary condition we can test in some cases. (This was
the approach taken in \cite{Milley}.)
\begin{theorem}[\cite{FKP} 1.1]
Suppose $M$ is a hyperbolic 3-manifold of finite volume.
For each component $k$ of $\partial \overline{M}$,
let $B_k$ be an embedded horoball neighborhood of $k$,
such that all the $B_k$ are disjoint.
Likewise, for each $k$,
let $s_k$ be a simple closed geodesic on $\partial B_k$.\footnote{With
$\partial B_k$ being the boundary of the completion
as above. However, this might in fact correspond to the
closure in $M$, as we have not assumed $B_k$ is maximal.}

Suppose the lengths of the $s_k$ are all greater than $2\pi$,
so that the Dehn filling $N$ of $M$ along all the $s_k$
is hyperbolic by Perelman and the $2\pi$-Theorem.

Let $\ell$ be the minimum among the lengths of all the $s_k$.

Then
\[
vol(N) \geq (1 - (2\pi / \ell)^2)^{3/2} \cdot vol(M).
\]
Equivalently,
\begin{equation}\label{eq:fkp}
\ell^2 \leq (2\pi)^2 / (1 - (vol(N)/vol(M))^{2/3}).
\end{equation}
\end{theorem}

This yields the following corollary.
\begin{corollary}
  Suppose $M,N$ are orientable, complete
  hyperbolic 3-manifolds of finite volume.
  If $vol(N)/vol(M) < 1$,
  then $N$ is a Dehn filling of $M$
  only if it is a Dehn filling of $M$ along a slope $s$
  whose minimum length $\ell$ satisfies equation (\ref{eq:fkp}).
\end{corollary}

When the conditions of this corollary are satisfied,
we can thus disprove $N$ being a Dehn filling of $M$
by enumerating all Dehn fillings of $M$ with small enough slope,
then for each Dehn filling, checking isometry with $N$.
If $N$ is isometric to none of these finitely many manifolds,
then $N$ is not a Dehn filling of $M$ at all.

On the other hand, if $vol(N)/vol(M) > 1$, 
then we know that $N$ is not a Dehn filling of $M$, because
\begin{theorem}[\cite{Th82}, 3.4 (c)]
Dehn filling decreases volume.
\end{theorem}

The upper bound on squared length in (\ref{eq:fkp})
is asymptotic to $(3/2)\cdot vol(N)/\Delta V$
as $\Delta V \to 0$. It would be of interest to get a bound independent
of $vol(N)$, and, especially, to get a left-hand bound as well. The
latter would reduce the number of cases needed to prove that
$N$ is not a Dehn filling of $M$.

In \cite{HK08} Hodgson and Kerckhoff gave a bound
that appears to completely solve the problem.
It uses the following notion\footnote{This
definition was anticipated in Theorem IA of \cite{NZ85}.
Hodgson and Kerckhoff use $\hat{L}$ instead
of $\mathbf{L}$, but we find $\hat{L}$ too
typographically close to $L$ for comfort.}
\begin{definition}
  The \emph{normalized length} $\mathbf{L}(s)$,
  of a single slope $s$
  on a horospherical cusp torus $k$ is
  \[
  \mathbf{L}(s) = \frac{length_k(s)}{\sqrt{area(k)}}.
  \]
  
  For a more general Dehn filling coefficient $c$,
  letting $c(k)$ be $c$'s slope choice on $k$,
  define $\mathbf{L}(c)$ by requiring $\mathbf{L} \geq 0$ and
  \[
  \frac{1}{\mathbf{L}(c)^2} = 
  \sum_{k \in K_c} \frac{1}{\mathbf{L}(c(k))^2}.
  \]
\end{definition}

Using this definition of normalized length,
they proved the following theorem:
\begin{theorem}[\cite{HK08}, Thm. 5.11, Cor. 5.13]\label{thm:mot}
  Let $M$ be a compact orientable 3-manifold
  whose interior admits a complete hyperbolic metric of finite volume.
  Let $c$ be a Dehn filling of $M$
  such that $\mathbf{L}(c) > 7.5832$. Then
  \begin{itemize}
  \item $M_c$ itself admits
    a complete hyperbolic metric on its interior;
  \item $M$ is homeomorphic to $M_c \smallsetminus \gamma$,
    where $\gamma$ is a geodesic link of $M_c$
    of total length at most $0.156012$
    admitting simultaneously embedded tubes
    about its components of radius at least
    $\artanh(1/\sqrt{3})$; and
  \item $\Delta V(c) < 0.198$.
  \end{itemize}
\end{theorem}

After suitably rephrasing this,
it gives a method to
test for Dehn filling---i.e.,
a Dehn ``parental test'':
\begin{corollary}
  Let $M$, $N$ be orientable 3-manifolds
  admitting complete hyperbolic metrics of finite volume
  on their interior.

  $N$ is a Dehn filling of $M$
  if and only if either
  \begin{itemize}
  \item $N$ is a Dehn filling of $M$
    along a slope of normalized length
    at most 7.5832, or
  \item $M$ is isometric to $N\setminus\gamma$
    for $\gamma$ a geodesic link
    of length less than 0.156012.
  \end{itemize}
\end{corollary}
The putative method runs as follows.
For simplicity's sake, let us assume
$M$ has exactly one more cusp than does $N$.
Suppose we wish to know if $N$ is a
Dehn filling of $M$, with $M,N$ both
finite-volume orientable hyperbolic 3-manifolds.
The collection $S$ of Dehn filling coefficients of $\partial M$
with normalized length less than 7.5832 is finite and computable,
since $M$ has exactly one more cusp than $N$ does.
The length spectrum of $N$ up to a given length $\ell$
is also a finite, computable set, as is the set $\Gamma_\ell$
of associated closed geodesics. Thus, to determine if $N$ is
a Dehn filling of $M$, it would suffice to check whether or not
$N$ is isometric to $M_s$ for all $s \in S$,
and to check whether or not $M$ is isometric
to $N \setminus \gamma$ for all $\gamma \in \Gamma_{0.156012}$.

The present work would be obviated by the above method
if there were efficient procedures to compute $\Gamma_\ell$
(or some reasonable superset thereof) for $\ell = 0.156012$,
to calculate the simple closed geodesics with given length,
and to calculate triangulations of $N \setminus \gamma$
for all simple $\gamma \in \Gamma_\ell$ (or the superset).

Unfortunately, there are no known, efficient such procedures.
The state of the art in hyperbolic 3-manifold software
is the remarkable program SnapPy \cite{SnapPy}.
SnapPy can only drill out simple closed curves
in the dual 1-complex of an ideal triangulation.
(See the file \texttt{drilling.c} of the SnapPy kernel
available at \cite{SnapPy}.) As explained in \cite{HW94}
on page 264, these dual 1-complex curves
may or may not be isotopic---or even homotopic---to
a given geodesic which one wishes to drill out.
Moreover, it is unclear not only
how to conveniently and efficiently
represent an arbitrary closed geodesic $\gamma$,
but also how to efficiently compute a triangulation
of $N \setminus \gamma$.
Furthermore, the current method
for calculating the length spectrum
(introduced in \cite{HW94})
depends crucially upon the construction
of a Dirichlet domain for an associated action
of a Kleinian group on hyperbolic space.
It is computationally expensive
to compute a reasonable guess
for a Dirichlet domain.
Finally, it is unclear
how to rigorously verify
the combinatorics and geometry of a
guess for the Dirichlet domain,
or to rigorously affirm or, especially, deny
the simplicity of a given closed geodesic,
without resorting to exact arithmetic,
which is significantly slower than an
approximate arithmetic like interval arithmetic.

The present work gets around these difficulties
by avoiding drilling altogether and by not requiring
a calculation of the length spectrum up to 0.156012
or even a calculation of the systole length, but only
requiring a lower bound on systole length (see
subsection \ref{sec:implement}).

\section{Rewriting the Hodgson-Kerckhoff Bounds}\label{sec:thm}

Hodgson and Kerckhoff's Theorem \ref{thm:mot}
is a corollary of their more powerful Theorem \ref{thm:hkbds}.
In this section we use Theorem \ref{thm:hkbds}
and other definitions and results from \cite{HK08}
to prove the main result of the paper.
We first review Theorem \ref{thm:hkbds},
which gives bounds on $\Delta V$ and $\ell$
in terms of auxiliary ``$z$''-variables.
In subsection \ref{sub:mon},
we show that the bounds
are decreasing functions of the $z$-variables. 
We also prove the monotonicities
of some auxiliary functions
necessary in what follows.
In subsection \ref{sub:complicated}
we rewrite the preceding bounds
into bounds on $\mathbf{L}(c)$ and $\ell$
in terms of $\Delta V$,
and not in terms of the $z$-variables.
These precise bounds are complicated,
so in subsection \ref{sub:nice} we relax them
to bounds which are more easy
to implement rigorously in code.
Finally, in subsection \ref{sub:num}
we conclude the proof of the main theorem.
Sample code that one could use to verify
these elementary considerations is available
at the author's GitHub repository \cite{bobbycyiii}.

Let us now restate Hodgson and Kerckhoff's theorem as we will use it.
\begin{theorem}[\cite{HK08}, 5.12]\label{thm:hkbds}
Let $M$ and $c$ be as in Theorem \ref{thm:mot}---in
particular, $\mathbf{L}(c) > 7.5832$.
Let $\Delta V = vol(M) - vol(M_c)$. Let $\ell$
be the total length of the geodesic link at the core of the filling.
Then
\begin{equation}\label{eqn:hkbds}
\frac{1}{4} \cdot \int_{ z_L}^1 
\frac{H'(z)}{H(z) \cdot (H(z) -  G_L(z))}\,dz
\leq \Delta V \leq
\frac{1}{4} \cdot \int_{ z_U}^1 
\frac{H'(z)}{H(z) \cdot (H(z) + G_U(z))}\,dz,
\end{equation}
and
\begin{equation}\label{eqn:lenbds}
1/H( z_L) \leq 2 \pi \cdot \ell
\leq 1/H( z_U),
\end{equation}
where $H,G, G_L, z_L$, and $ z_U$
have the following definitions.
\end{theorem}
\begin{definition}\label{def:vars}
\begin{align*}
S = \frac{\frac{1}{2\cdot\sqrt{2}}}{\arsinh \left( \frac{1}{2\cdot\sqrt{2}} \right) },&\qquad
K = \frac{2\sqrt{3}}{S} \approx 3.3957,\\
h(z) = \frac{1+z^2}{z\cdot (1 - z^2)},\qquad
g_U(z) &= \frac{1+z^2}{2\cdot z^3},\qquad
g_L(z) = \frac{(1+z^2)^2}{2\cdot z^3 \cdot (3 - z^2)},\qquad\\
H = h/K,\qquad G_U &= g_U/K, \qquad G_L = g_L/K,
\end{align*}
\begin{align*}
F_U(z) = \frac{H'(z)}{H(z) + G_U(z)} - \frac{1}{1-z}
&= \frac{h'(z)}{h(z) + g_U(z)} - \frac{1}{1-z},\\
 F_L(z) = \frac{H'(z)}
                    {H(z) -  G_L(z)} - \frac{1}{1-z}
&= \frac{h'(z)}{h(z) - g_L(z)} - \frac{1}{1-z},\\
f_U(z) = K \cdot (1-z) \cdot e^{-\Phi_U(z)},&\quad
\Phi_U(z) = \int_1^z F_U(w)\,dw,\\
 f_L(z)=K \cdot (1-z) \cdot e^{- \Phi_L(z)},&\quad
 \Phi_L(z) = \int_1^z  F_L(w)\,dw,\\
f_U( z_U) = \frac{(2\pi)^2}{\mathbf{L}(c)^2}, \qquad
 &f_L( z_L) = \frac{(2\pi)^2}{\mathbf{L}(c)^2}.
\end{align*}
\end{definition}
\begin{remark}
These definitions are from 
pp. 1079, 1080, 1084, and 1088 of \cite{HK08}. 
We have taken the liberty of making the following
changes to the notation, eliminating 
the small over-accents on the letters:
\[
\begin{matrix}
f_U \gets f, & f_L \gets \tilde{f}, & F_U \gets F, & F_L \gets \tilde{F},\\
g_U \gets g, & g_L \gets \tilde{g}, & G_U \gets G, & G_L \gets \tilde{G},\\
\Phi_U \gets \Phi, & \Phi_L \gets \tilde{\Phi}, & z_L \gets \tilde{z}, & z_U \gets \hat{z}.
\end{matrix}
\]
The reader should note that the above theorem has
$2 \pi \cdot \ell$ in place of $\mathcal{A}$.
This is valid---see, e.g., Corollary 5.13 of \cite{HK08}.
\end{remark}

\begin{remark}
The considerations in \cite{HK08}\ 
in pp. 1079--1081 show that, under the
assumption $\mathbf{L}(c) > 7.5832$,
we have $z_U,z_L \in (\sqrt{1/3},1)$.
Hodgson and Kerckhoff's upper bound
$0.198$ on $\Delta V$ under the same
assumption is an upper bound on
the right-hand integral in equation
\ref{eqn:hkbds} with $z_U$ replaced with
$\sqrt{1/3}$, i.e. $UB(\sqrt{1/3})$ with
the definition below.
\end{remark}

\subsection{Monotonicities}\label{sub:mon}
Let
\begin{equation}\label{eqn:lb}
LB(z) = \frac{1}{4} \cdot \int_{z}^1 
\frac{H'(w)}{H(w) \cdot (H(w) -  G_L(w))}\,dw
\end{equation}
and
\begin{equation}\label{eqn:ub}
UB(z) = \frac{1}{4} \cdot \int_{z}^1 
\frac{H'(w)}{H(w) \cdot (H(w) + G_U(w))}\,dw.
\end{equation}

We intend to transform the bounds on
$LB$ and $UB$ given in Theorem \ref{thm:hkbds}
into bounds on $z_L$ and $z_U$
in terms of $\Delta V$---and
thence into bounds on $\mathbf{L}(c)$ and $\ell$---by 
inverting $LB$ and $UB$.
This will work if
we know the monotonicity of $LB$ and $UB$.
We will also later require the monotonicities
of $H$, $f_L$, and $f_U$,
so we also prove their monotonicities in this section.
The proofs of the following lemmas
are all straightforward, similar calculations.
We give the longest proof in detail,
and give more abbreviated proofs for the other lemmas.
We omit the proofs of Lemmas \ref{lem:f}\ and \ref{lem:ftilde},
as they are proved in \cite{AtkFut}.

\begin{lemma}\label{lem:lb}
$LB$ is decreasing on $\left(\sqrt{\sqrt{5}-2},1\right)$.
\end{lemma}

\begin{proof}
  By the definition of $LB$, we have
  \[
  LB'(z) = 
  -\frac{1}{4} \cdot 
  \frac{H'(z)}{H(z)\cdot(H(z)- G_L(z))}.
  \]
  This is equal to
  \[
  - \frac{K}{4} \cdot \frac{h'(z)}{h(z)\cdot(h(z)-g_L(z))}.
  \]
  Since $K > 0$, this has the same sign as
  \[
  - \frac{h'(z)}{h(z)\cdot(h(z)-g_L(z))}.
  \]
  Expanding the definitions of $g$ and $h$ and factoring,
  this is equal to
  \[
  -\frac{2\cdot z^2\cdot(z^2-3)\cdot(z^4+4\cdot z^2 -1)}
  {(z^2+1)^2\cdot(z^2-2\cdot z-1)\cdot(z^2+2\cdot z-1)}.
  \]
  Since $|z| < 1$, we have $z^2 - 3 < 0$.
  Therefore, the above has the same sign as
  \[
  \frac{z^2\cdot(z^4+4\cdot z^2 -1)}
       {(z^2+1)^2\cdot(z^2-2\cdot z-1)\cdot(z^2+2\cdot z-1)}.
  \]
  Since $z > \sqrt{\sqrt{5}-2}$, we also have $z^2/(z^2+1)^2 > 0$.
  Therefore, the above has the same sign as
  \[
  \frac{z^4+4\cdot z^2 -1}
       {(z^2-2\cdot z-1)\cdot(z^2+2\cdot z-1)}.
  \]
  For the same reason, $z^4+4\cdot z^2-1>0$,
  and therefore the above has the same sign as
  \[
  \frac{1}{(z^2-2\cdot z-1)\cdot(z^2+2\cdot z-1)}.
  \]
  Since $\sqrt{\sqrt{5}-2} > \sqrt{2} - 1$,
  we have $z^2 + 2\cdot z - 1 > 0$;
  therefore the above has the same sign as
  \[
  \frac{1}{z^2-2\cdot z-1}.
  \]
  Finally, $z > 1-\sqrt{2}$ because $z > \sqrt{\sqrt{5}-2}$,
  and  $z < 1 + \sqrt{2}$ because $z < 1$;
  therefore $z^2-2\cdot z - 1 < 0$,
  so that the above is negative.
  Since $LB'$ is negative on $\left(\sqrt{\sqrt{5}-2},1\right)$,
  $LB$ is decreasing on this interval.
\end{proof}

\begin{lemma}\label{lem:ub}
$UB$ is decreasing on 
$\left(\sqrt{\sqrt{5}-2},\infty\right)$.
\end{lemma}
\begin{proof}
Simplifying $UB'(z)$ yields the expression
\[
-\frac{K}{2}\cdot
\frac{z^2 \cdot (z^4 + 4\cdot z^2 - 1)}{(z^2 + 1)^3};
\]
on the given domain, $z^4 + 4\cdot z^2 - 1$ is positive,
and hence $UB'(z) < 0$.
\end{proof}

\begin{lemma}\label{lem:H}
$H$ is increasing on 
$\left(\sqrt{\sqrt{5}-2},\infty\right)$.
\end{lemma}
\begin{proof}
$H'(z)$, after multiplying by $K > 0$, simplifies to
\[
\frac{z^4 + 4\cdot z - 1}
     {(z-1)^2\cdot z^2 \cdot (z+1)^2};
\]
since $z^4 + 4\cdot z - 1$ is positive on the
given domain, $H'(z) > 0$ there.
\end{proof}

\begin{lemma}\label{lem:f}
$f_L$ is decreasing on 
$\left(\sqrt{\sqrt{5}-2},\infty\right)$.
\end{lemma}
\begin{lemma}\label{lem:ftilde}
$f_U$ is decreasing on 
$\left(\sqrt{\sqrt{5}-2},\sqrt{3}\right)$.
\end{lemma}
\begin{proof}[Proof of Lemmas \ref{lem:f} and \ref{lem:ftilde}.]
See \cite{AtkFut}, Lemma 5.4.
\end{proof}

Thus $LB,$ $UB,$ $f_L,$ $f_U,$ and $H$
are invertible on the given intervals.
Since $[\sqrt{1/3},1]$ is a subset
of each of the given intervals,
we may invert these functions on $[\sqrt{1/3},1]$.
This is important because
the arguments to these functions are $z$-variables,
which in \cite{HK08} (see p. 1079) are assumed
always to lie in $[\sqrt{1/3},1].$

\begin{definition}
Let $BU$ be the inverse of the restriction
of $UB$ to the interval $[\sqrt{1/3},1]$.
Likewise, let $BL$ be the inverse of the
restriction of $LB$ to the interval
$[\sqrt{1/3},1]$.
\end{definition}

\subsection{Explicit bounds on $\mathbf{L}(c)$ 
            and $\ell$ in terms of $\Delta V$}\label{sub:complicated}

Having established these monotonicity results,
let us now derive bounds on $\mathbf{L}(c)$ and $\ell$
in terms of $\Delta V$.

We first require the following surprising lemma.
\begin{lemma}\label{lem:LBgtUB}
  For $z \in (\sqrt{1/3},1)$,
  $LB(z) > UB(z)$.
\end{lemma}
This is not in contradiction with
Theorem \ref{thm:hkbds},
which states $LB(z_L) \leq \Delta V \leq UB(z_U)$
for particular $z_L, z_U$, and
not $LB(z) \leq \Delta V \leq UB(z)$.
\begin{proof}
  The difference of the integrands in the definitions
  of $LB$ and $UB$, in that order, simplifies to
  \[
  8\cdot K \cdot
  \frac{(z-1)\cdot z^2 \cdot (z+1) \cdot (z^4 + 4\cdot z^2 - 1)}
       {(z^2+1)^3 \cdot (z^2 - 2\cdot z - 1) \cdot (z^2 + 2\cdot z -1)}.
  \]
  By the proofs of the monotonicities above,
  all factors except for $z-1$ and $z^2-2\cdot z -1$
  in this expression are positive for $z \in (\sqrt{1/3},1)$.
  Therefore this expression is positive on this interval.
  So the integrand of $LB$ is greater than the
  integrand of $UB$ on this interval.
  Thus by the monotonicity of integration,
  $LB(z) > UB(z)$ for $z \in (\sqrt{1/3},1)$.  
\end{proof}
Now we can derive the explicit bounds desired.
\begin{lemma}\label{lem:nab}
Suppose $M$ and $c$ are as in Theorem \ref{thm:mot}.
Abusing notation, write $\Delta V$ for $\Delta V(c)$
and $\mathbf{L}$ for $\mathbf{L}(c)$. Then
\begin{equation}\label{eqn:Lbds}
\frac{(2\pi)^2}{f_L(BL(\Delta V))}
\leq
\mathbf{L}^2
\leq
\frac{(2\pi)^2}{f_U(BU(\Delta V))}.
\end{equation}
\end{lemma}
\begin{proof}
We know by definition that
$\frac{(2\pi)^2}{\mathbf{L}^2} = f_U(z_U) =  f_L(z_L)$.

To get an upper bound on $\mathbf{L}$,
we can get a lower bound on $f_U(z_U)$,
which would result from an upper bound on $z_U$
(since $f_U$ is decreasing),
which itself would result from a lower bound on $UB(z_U)$
(since $UB$ is decreasing).
Now, $\Delta V \leq UB( z_U)$ by Theorem \ref{thm:hkbds}.
Since $\mathbf{L}(c) > 7.5832$,
$\Delta V < UB(\sqrt{1/3})$.
So we can invert, and $BU(\Delta V) \geq z_U$,
since $BU$ is decreasing.
Then $f_U(BU(\Delta V)) \leq f_U(z_U) = (2\pi)^2 / \mathbf{L}^2$,
since $f_U$ is also decreasing. Therefore,
\[
\mathbf{L}^2 \leq
\frac{(2\pi)^2}
     {f_U(BU(\Delta V))},
\]
as desired.
One derives the lower bound similarly
from $LB(z_L) \leq \Delta V$.
We can invert since $\Delta V < UB(\sqrt{1/3}) < LB(\sqrt{1/3})$,
the last inequality coming from lemma \ref{lem:LBgtUB}.
\end{proof}

We conclude this subsection with the bound for $\ell$.
By equation \ref{eqn:lenbds} of Theorem \ref{thm:hkbds},
we already have an upper bound on $\ell$,
viz. $\ell \leq 1/(2\pi \cdot H( z_U))$.
We just need to put the right-hand side
in terms of $\Delta V$.

In fact, since $H$ is increasing,
$1/(2\pi \cdot H)$ is decreasing.
Therefore we just need a lower bound on $ z_U$;
applying $1/(2\pi \cdot H)$ to this lower bound
will give us an upper bound on $\ell$.

At this point, one could use the standing assumption
that $z$-variables have $z \in [\sqrt{1/3},1]$.
As a matter of fact, this is where the bounds in
Theorem \ref{thm:mot} come from.
But we would like a better bound for small $\Delta V$.
So in the next lemma we assume a
stronger hypothesis on $\Delta V$, which we ensure later
in our Main Theorem by assuming a hypothesis
on the length of the Dehn filling coefficient that
is stronger than what Hodgson and Kerckhoff assume.
Happily, the resulting bounds are on the same order of magnitude
as Hodgson and Kerckhoff's bounds.

\begin{lemma}
  Suppose $M$ and $c$ are as in Theorem \ref{thm:mot}, with
  $\ell$ the total length of the geodesic link at the core of the filling $M_c$.
  Abusing notation as in lemma \ref{lem:nab}, and assuming 
  $f_L (BL (\Delta V)) < f_U (\sqrt{1/3})$, we have
  \begin{equation}\label{eqn:ellbd}
    \ell \leq
    \frac{1}{2\pi\cdot (H \circ f_U^{-1} \circ  f_L \circ BL)
                       (\Delta V)},
  \end{equation}
  where $f_U^{-1}$ is taken to have domain $(0, f_U(\sqrt{1/3}))$.
\end{lemma}
\begin{proof}
We know $f_U(z_U) =  f_L(z_L)$ from Definition \ref{def:vars}.
Now, $f_L$ and $f_U$ both are decreasing. So
if we can get a lower bound on $z_L$, we can get
a lower bound on $z_U$, via upper bounds
on $f_U(z_U) =  f_L(z_L)$.
The left inequality of equation (\ref{eqn:hkbds})
from Theorem \ref{thm:hkbds}
says $LB(z_L) \leq \Delta V$; since $LB$ is 
decreasing, we get $z_L \geq BL(\Delta V)$.
As desired, this yields the lower bound
$z_U \geq f_U^{-1} (f_L (BL (\Delta V)))$.
Thus, by equation \ref{eqn:lenbds} of Theorem \ref{thm:hkbds},
since $H$ is increasing,
\[
\ell \leq \frac{1}{2\pi\cdot H(z_U)} \leq \frac{1}
{2\pi\cdot (H \circ f_U^{-1} \circ  f_L \circ BL)
           (\Delta V)},
\]
assuming it is well-defined, i.e. assuming
$f_L(BL(\Delta V))$ is in the domain of $f_U$.
This is precisely the last assumption of the lemma.
\end{proof}

\begin{remark}\label{rmk:theta}
The assumption $f_L(BL(\Delta V)) < f_U(\sqrt{1/3})$
is satisfied if 
$\Delta V < \Theta = LB(f_L^{-1}(f_U(\sqrt{1/3})))$,
which we can achieve by assuming
$L(c) > 2\cdot \pi / \sqrt{f_U(BU(\Theta))}$.
One may estimate $\Theta < 0.1562$ and
$2\cdot \pi / \sqrt{f_U(BU(\Theta))} < 8.5$.
\end{remark}
The reader may find better bounds
using the numerical estimation methods
of a computer algebra system. However, the
rigor of such computations is open to doubt.
The above are bounds we could prove rigorously
in the proof assistant Coq \cite{Coq} using its
interval arithmetic module \cite{CoqInterval}.
Verifications of these bounds are available
at the author's GitHub repository \cite{bobbycyiii}.
  
\subsection{Nice bounds}\label{sub:nice}
One could, conceivably, use equations
\ref{eqn:Lbds}\ and\ \ref{eqn:ellbd} to
get a Dehn parental test. However, the
expressions in the test are so complicated
that running the test in floating-point
arithmetic likely would yield meaningless results,
and running the test in interval arithmetic would not
reach sufficient precision in a reasonable
amount of time. So we relax these bounds
into expressions more quickly and precisely
computable. We first relax the bounds to
expressions with a less complicated dependence
on $z$. Then we finally derive appropriate numerical
approximations to the constants in the expressions.

The conditions which the approximations
should satisfy are not difficult to
derive. For instance, an approximation
$\eta$ to $1/(2\pi \cdot H)$ should
be decreasing, since $1/(2\pi\cdot H)$
is itself decreasing and we want a
reasonable approximation; and $\eta$
should be at least $1/(2\pi\cdot H)$
so that we can deduce 
\[
\ell \leq
(\eta \circ f_U^{-1} \circ  f_L \circ BL)(\Delta V)
\]
from (\ref{eqn:ellbd}). In fact, 
$\eta(z) = K\cdot(1-z)/(2\pi)$
suffices. As it turns out, we can approximate
all the necessary functions by linear functions of the form
$z \mapsto \kappa \cdot (1-z)$ for some $\kappa > 0$.
Those approximations are as follows:

\begin{lemma}\label{lem:bounds}
Let $Z = [\xi,1]$, where $\xi = BL(\Theta) \in [0.8112, 0.8113]$. 
Then for all $z \in Z$,
\begin{align}
1/h(z) &\leq 1-z,\label{eqn:rechbd}\\
f_U(z) &\geq A \cdot (1-z),\label{eqn:fUbd}\\
f_L(z) &\leq B\cdot(1-z),\label{eqn:fLbd}\\
LB(z) &\geq C\cdot(1-z),\label{eqn:LBbd}\\
UB(z) &\leq D\cdot(1-z),\label{eqn:UBbd}
\end{align}
where
\[
A = \frac{f_U(\xi)}{1 - \xi}\mbox{, }
B = \frac{f_L(\xi)}{1 - \xi}\mbox{, }
C = \frac{\Theta}{1 - \xi}\mbox{, and }
D = K/4.
\]
\end{lemma}

We thank the referee for suggesting the use
of secant line approximations instead
of Taylor approximations. Where they work,
they afford better and easier approximations.

\begin{proof}[Proof of \ref{eqn:rechbd}.]
  \[
  1 - z - \frac{1}{h(z)} = \frac{(1-z)^2}{1+z^2} \geq 0. \qedhere
  \]
\end{proof}

\begin{proof}[Proof of \ref{eqn:fUbd}.]
  Recall that
  \[
  f_U(z) = K \cdot (1-z) \cdot e^{-\Phi_U(z)}.
  \]
  Therefore,
  \[
  f_U''(z) = K \cdot e^{-\Phi_U(z)} \cdot
  \left (
  (1-z) \cdot \left ( \Phi_U'(z)^2 - \Phi_U''(z) \right ) + 2 \cdot \Phi_U'(z) \right).
  \]
  Now, $\Phi_U' = F_U$. Thus the third factor is
  \[
  \frac{f_U''(z)}{K \cdot e^{-\Phi_U(z)}} =
  (1-z) \cdot \left ( F_U(z)^2 - F_U'(z) \right ) + 2 \cdot F_U(z).
  \]
  Calculating, this is equal to
  \begin{equation}\label{expr:fUf3}
  -2 \cdot \frac{z^8+6\cdot z^6 + 32\cdot z^4 + 10\cdot z^2 - 1}
                {(z+1)\cdot (z^2 + 1)^4}.
  \end{equation}
  The denominator is clearly positive on $Z$. The numerator
  has two real roots, neither of which lies in $Z$. (One can
  estimate these roots by first isolating them in intervals,
  e.g. via the slow but easy method of Sturm, then refining
  the intervals, e.g. by bisection.) So it
  has constant sign on $Z$. By inspection at 1, that sign is positive.
  Therefore $f_U'' < 0$ on $Z$, so $f_U$ is concave down. So the
  secant line over $Z$ is a lower bound for $f_U$ on $Z$.
\end{proof}

\begin{proof}[Proof of \ref{eqn:fLbd}.]
  As above in the proof of \ref{eqn:fUbd}, 
  the concavity of $f_L$ reduces to determining
  the sign of
  \[
  \frac{f_L''(z)}{K \cdot e^{-\Phi_L(z)}} =
  (1-z) \cdot \left ( F_L(z)^2 - F_L'(z) \right ) + 2 \cdot F_L(z).
  \]
  Calculating, this is equal to
  \begin{equation}\label{expr:fLf3}
  -2 \cdot 
  \frac{z^{12} -4\cdot z^{10} +17\cdot z^8 -248\cdot z^6 +203\cdot z^4 -36 \cdot z^2 +3}
       {(z+1)\cdot(z^2+1)^2\cdot(z^2-2\cdot z-1)^2\cdot(z^2+2\cdot z-1)^2}.
  \end{equation}
  The denominator is not zero on $Z$, and therefore
  is clearly positive on $Z$. Likewise, as above,
  after estimating the roots of the numerator,
  we can prove none lies in $Z$. So expression
  \ref{expr:fLf3} has constant sign on $Z$.
  Surprisingly, it evaluates to 1 at $z=1$. Therefore,
  on $Z$, $f_L$ is concave up, and the secant line over
  $Z$ is an upper bound to $f_L$ over $Z$.
  
  The secant line approximation has slope
  $LB(\xi)/(1-\xi)$, but $LB(\xi) = \Theta$.
  \end{proof}
  
\begin{proof}[Proof of \ref{eqn:LBbd}.]
  We aim to show that the secant line over $Z$ suffices
  by showing $LB$ is concave down over $Z$.
  
  Let
  \[
  lb(z) = \int_z^1 \frac{h'(w)}{h(w)\cdot(h(w)-g_L(w))}\,dw.
  \]
  Then $K\cdot lb(z) / 4 = LB(z).$ It will thus suffice 
  to show $lb$ is concave down.

  Let $t_L$ be the integrand $h'/(h\cdot (h - g_L))$.
  Then $lb''(z) = -t_L'(z)$. Factoring, we get
  \[
  t_L'(z) = - \frac{4 \cdot (z-1) \cdot z \cdot (z+1) \cdot p(z)}
  {(z^2+1)^3\cdot(z^2-2\cdot z-1)^2\cdot
  (z^2+2\cdot z-1)^2},
  \]
  where
  \[
  p(z) = 5 \cdot z^8 - 6 \cdot z^6 + 88 \cdot z^4
  - 26\cdot z^2 + 3.
  \]
  Now, $t_L'$ has the same sign as $p$ on the
  interior $Z'$ of $Z$, since $-1$ and $z-1$ are negative
  on this interior, and all the other factors are
  positive on $Z$. One may calculate that $p$ has
  no real roots whatever. For variety's sake, we
  do this the following way.
  \begin{align*}
    p(z)
    & =
    5 \cdot z^8 - 6 \cdot z^6 + 2 \cdot z^4
    + 86 \cdot z^4 - 26 \cdot z^2 + 3 \\
    & =
    z^4\cdot(5 \cdot (z^2)^2 - 6 \cdot (z^2) + 2)
    +
    86 \cdot (z^2)^2 - 26 \cdot z^2 + 3.
  \end{align*}
  Now, $(-6)^2 - 4\cdot 5 \cdot 2 < 0$ and
  $(-26)^2 -4 \cdot 86 \cdot 3 < 0$.
  Therefore, $5\cdot z^2 - 6\cdot z +2$ has constant sign,
  and $86 \cdot z^2 - 26 \cdot z + 3$ does too. By evaluation
  at 0, this sign is positive on both. Therefore $p$ is positive.
  (Alternatively, one can again use a root isolation
  algorithm as in the previous lemmas.) 
  In any case, $t_L' > 0$ on $Z$. So $lb''(z) < 0$ on $Z$.
  Thus, finally, $LB$ is concave down on $Z$, since $LB$ is a
  positive constant multiple of $lb$.
\end{proof}

\begin{proof}[Proof of \ref{eqn:UBbd}.]
  We would like to prove the secant line approximation suffices
  by showing $UB$ is concave up on $Z$. As above,
  we define
  \[
  ub(z) = \int_z^1 \frac{h'(w)}{h(w)\cdot(h(w)+g_U(w))}\,dw,
  \]
  so that $K\cdot ub(z)/4 = UB(z)$. It would suffice
  to show $ub$ is concave up.

  Let $t_U$ be the integrand $h'/(h\cdot (h + g_U)).$ Then
  $ub''(z) = -t_U'(z)$. Factoring, we get
  \[
  t_U'(z) =
  - \frac
  {4 \cdot z \cdot (z^4 - 10 \cdot z^2 + 1)}
  {(z^2 + 1)^4}.
  \]
  We see that $t_U'$ has the opposite sign of
  $p(z) = z^4 - 10 \cdot z^2 + 1$
  on $Z$. None of the four real roots of $p$
  lies in $Z$, and so $p$ has constant sign
  on $Z$. By inspection at $1$, this sign is
  negative. Therefore, $t_U' > 0$ on $Z$. But
  that means $ub'' < 0$ on $Z$, and in fact
  $ub$ is concave \emph{down} on $Z$, contrary
  to what we would have liked.

  However, since $ub$ is concave down, we can
  use the linearization at 1 as an upper bound.
  The slope is, surprisingly, $-t_U(1) = -1$ at 1.
  So $1-z$ is an upper bound to $ub$ on $Z$,
  and hence $(K/4)\cdot (1-z)$ is an upper bound
  to $UB$ on $Z$.
\end{proof}

This leaves us with the following nice bounds.

\begin{lemma}\label{lem:nicebds}
  If $0 < \Delta V < \Theta$, then
\begin{equation}\label{eqn:lnice}
\frac{1}{2 \pi \cdot (H \circ f_U^{-1} \circ f_L \circ BL)(\Delta V)}
\leq \alpha \cdot \Delta V,
\end{equation}
\begin{equation}\label{eqn:lclnice}
\frac{(2 \pi)^2}{ f_L(BL(\Delta V))}
\geq \beta \cdot \frac{1}{\Delta V},
\end{equation}
and
\begin{equation}\label{eqn:lcunice}
\frac{(2 \pi)^2}{f_U(BU(\Delta V))}
\leq \gamma \cdot \frac{1}{\Delta V},
\end{equation}
where, letting $\xi = BL(\Theta),$
\[
\alpha = \frac{K \cdot (1-\xi)}{2\cdot \pi \cdot \Theta}
  \cdot e^{\Phi_U(\xi) - \Phi_L(\xi)},
\qquad
\beta = \frac{(2 \pi)^2\cdot \Theta}
             {f_L(\xi)},
\qquad
\mbox{and }
\gamma = \pi^2 \cdot e^{\Phi_U(\xi)}.
\] 

\end{lemma}
\begin{proof}[Proof of \ref{eqn:lnice}.]
  Recall that $H = h/K$. Thus the left-hand
  side of equation \ref{eqn:lnice} is
  \[(K/(2\cdot \pi)) \cdot
  ((1/h)\circ f_U^{-1} \circ f_L \circ BL) (\Delta V). \]
  By equation \ref{eqn:LBbd}, $LB(z) \geq C \cdot (1-z)$,
  so $BL(\Delta V) \geq 1 - \Delta V/C$. By equation \ref{eqn:fLbd},
  since $f_L$ is decreasing,
  \begin{align*}
    f_L(BL(\Delta V))
    & \leq f_L(1-\Delta V/C)\\
    & \leq B\cdot(1-(1-\Delta V/C))\\
    & = B\cdot \Delta V/C.
    \end{align*}
  By equation \ref{eqn:fUbd}, $f_U(z) \geq A \cdot (1-z),$
  so $f_U^{-1}(\Delta V) \geq 1 - \Delta V/A$.
  Since $f_U$ is decreasing,
  so is $f_U^{-1}$, and therefore
  \begin{align*}
    f_U^{-1}(f_L(BL(\Delta V)))
    & \geq f_U^{-1}(B\cdot \Delta V/C)\\
    & \geq 1-(B\cdot \Delta V/C)/A\\
    & = 1-B\cdot \Delta V/(A\cdot C).
    \end{align*}
  Since $H$ is increasing, $1/h$ is decreasing.
  By equation \ref{eqn:rechbd},
  \begin{align*}
    1/h(f_U^{-1}(f_L(BL(\Delta V))))
    & \leq 1/h(1-B\cdot \Delta V/(A\cdot C))\\
    &\leq 1-(1-B\cdot \Delta V/(A\cdot C))\\
    &= B\cdot \Delta V/(A \cdot C).
    \end{align*}
  Finally, since $K/(2\cdot \pi) > 0$, we have
  that the left-hand side of equation \ref{eqn:lnice}
  is at most
  $K \cdot B \cdot \Delta V/ (2\cdot \pi \cdot A \cdot C).$
  Calculating,
  \[\frac{K \cdot B}{2\cdot \pi \cdot A \cdot C}
  = \frac{K \cdot f_L(\xi) \cdot (1 - \xi)}
  {2\cdot \pi \cdot f_U(\xi) \cdot \Theta}
  = \frac{K \cdot e^{-\Phi_L(\xi)} \cdot (1 - \xi)}
  {2\cdot \pi \cdot e^{-\Phi_U(\xi)} \cdot \Theta}
  =\frac{K \cdot (1-\xi)}{2\cdot \pi \cdot \Theta}
  \cdot e^{\Phi_U(\xi) - \Phi_L(\xi)}. \qedhere \]
\end{proof}

\begin{proof}[Proof of \ref{eqn:lclnice}]
  As above, $BL(\Delta V) \geq 1 - \Delta V / C$.
  Since $f_L$ is decreasing,
  $f_L(BL(\Delta V))$ is at most $f_L(1-\Delta V/C)$,
  which, by equation \ref{eqn:fLbd}, is at most
  $B\cdot(1-(1-\Delta V/C)) = B\cdot \Delta V / C.$
  Thus $(2\cdot\pi)^2/f_L(BL(\Delta V))$ is at least
  $(2\cdot \pi)^2 \cdot C / (B \cdot \Delta V)$.
  But $C / B = \Theta / f_L(\xi)$.
\end{proof}

\begin{proof}[Proof of \ref{eqn:lcunice}]
  By equation \ref{eqn:UBbd}, as above,
  $BU(\Delta V) \leq 1 - \Delta V / D$.
  Since $f_U$ is decreasing,
  $f_U(BU(\Delta V))\geq f_U(1-\Delta V/D)$, which, by
  equation \ref{eqn:fUbd}, is, as above,
  at least $A\cdot \Delta V/D.$ Thus
  $(2\cdot \pi)^2/f_U(BU(\Delta V)) \leq (2\cdot \pi)^2 \cdot D / (A \cdot \Delta V).$
  But $D/A = K\cdot(1-\xi)/(4\cdot f_U(\xi))$, which
  is $e^{\Phi_U(\xi)}/4.$
\end{proof}

\subsection{Numerical approximations}\label{sub:num}
The bounds from Lemma \ref{lem:nicebds} could be
implemented as they stand,
but any software implementation thereof might spend
a large amount of time working out estimates of the
constant factors $\alpha$, $\beta$, and $\gamma$.
So we relax the bounds by getting simple, once-for-all
estimates on these factors.
Using Coq one may show, following
the code at \cite{bobbycyiii}, that
\begin{lemma}\label{lem:appx}
$\alpha < 0.97$,
$\beta > 8.95,$ and 
$\gamma < 13.4$.
\end{lemma}

\subsection{Proof of main theorem}\label{sub:main}

This is enough to complete the proof of the main theorem
of this paper. That theorem runs as follows.

\begin{main}
Let $M,N$ be orientable 3-manifolds
admitting complete hyperbolic metrics
of finite volume on their interiors.
Let $\Delta V = Vol(M) - Vol(N)$.

$N$ is a Dehn filling of $M$ if and
only if either
\begin{itemize}
\item $N$ is a Dehn filling of $M$
along a coefficient of normalized length
at most $8.5$, or
\item $N$ both \emph{has} a closed simple geodesic
of length less than $0.97 \cdot \Delta V$, and
\emph{is} a Dehn filling of $M$ along
a coefficient of normalized length $\mathbf{L}$
such that
\begin{equation}\label{eqn:finalbds}
\frac{8.95}{\Delta V} < \mathbf{L}^2 < \frac{13.4}{\Delta V}.
\end{equation}
\end{itemize}
\end{main}

\begin{proof}[Proof of Main Theorem.]
The if-direction is plain. For the only-if direction,
suppose $N$ is a Dehn filling of $M$. Then
either $N$ is a Dehn filling of $M$ along a
slope $c$ with $\mathbf{L}(c) \leq 8.5$
or $N$ is Dehn filling of $M$ along a slope $c$
with $\mathbf{L}(c) > 8.5$. The former case
is the first disjunct in the theorem.
In the latter case, Theorem \ref{thm:mot} applies, and
$\Delta V < \Theta < 0.1562$ by Remark \ref{rmk:theta}.
We can therefore apply Lemma \ref{lem:nicebds}. 
Now, by equations \ref{eqn:ellbd} and 
\ref{eqn:lnice}, the geodesic link at the core
of the filling has length $\ell$ satisfying
$\ell < \alpha \cdot \Delta V < 0.97 \cdot \Delta V$. 
So $N$ has a closed simple geodesic of length less than $0.97 \cdot \Delta V$.
Furthermore, by equations \ref{eqn:Lbds} and \ref{eqn:lclnice},
$8.95 / \Delta V < \beta / \Delta V \leq \mathbf{L}(c)^2$,
and by equations \ref{eqn:Lbds} and \ref{eqn:lcunice}, 
$\mathbf{L}(c)^2 \leq \gamma / \Delta V < 13.4 /\Delta V$.
So $N$ is in fact a Dehn filling of $M$ along a coefficient
of normalized length satisfying equation \ref{eqn:finalbds}.
Since $N$ has an appropriately short geodesic and is an
appropriately short Dehn filling of $M$, $N$ satisfies
the last disjunct of the theorem.
\end{proof}

\section{Prospects}\label{sec:prospects}
\subsection{Implementation}\label{sec:implement}
For rigorous estimates on numerical properties of hyperbolic
3-manifolds, formerly one had to compute using exact
arithmetic in the software suite Snap (see \cite{Snap}\footnote{Snap
  is no longer hosted at the University of Melbourne, but
  instead at \url{http://snap-pari.sourceforge.net}}),
and then get bounds on the exact computed numbers.
The state of the art today is to compute estimates
using an approximate arithmetic like interval arithmetic,
as in the HIKMOT module of SnapPy (see \cite{HIKMOT}, \cite{SnapPy}).
The main theorem of this paper lays the foundations
for implementing a Dehn parental test using
such an approximate arithmetic. But such an
algorithm does not immediately follow from these bounds.
The following is a sketch of an algorithm that works
under a simplifying assumption.

\subsubsection*{Parental test sketch with one cusp difference}
Suppose $M$ has one more cusp than $N$, and both are orientable
cusped hyperbolic 3-manifolds of finite volume.
\begin{enumerate}
\item Let $\Delta V = vol(M) - vol(N)$.
  Let $X$ be a finite superset of the Dehn
  filling coefficients of $M$ with
  normalized length at most $8.5$. (We can calculate
  such an $X$ using HIKMOT.)
\item For every $c \in X$, if $N$ is isometric to $M(c)$,
  then $N$ is a Dehn filling of $M$.
\item Otherwise, $N$ is not a filling along a short slope.
  So $N$ must have a short geodesic. Calculate a positive lower
  bound $S$ on the systole length of $N$.
  In interval arithmetic we can verify either
  $S > 0.97 \cdot \Delta V$ or $S < 0.98 \cdot \Delta V$.
  \begin{itemize}
  \item[($>$)] If $S > 0.97\cdot \Delta V$, then $N$
    is not a Dehn filling of $M$.
  \item[($<$)] If $S < 0.98\cdot \Delta V$, then $\Delta V > S/0.98 > 0$.
    The set $Y$ of all Dehn filling coefficients $c$
    with normalized length satisfying \ref{eqn:finalbds}
    is well-defined. Moreover, since $M$ has exactly one
    more cusp than $N$, $Y$ is finite. We can construct
    a finite superset $Y'$ using HIKMOT.
    Iterating over $Y'$, determine if
    there exists $c \in Y'$ such that $N$ is isometric to $M(c)$.
    If so, $N$ is a Dehn filling of $M$.
    Otherwise, $N$ is not a Dehn filling of $M$.
  \end{itemize}
\end{enumerate}

\subsubsection*{Discussion}
The first difficulty that a naive implementation
would encounter is that if the putative
parent manifold $M$ has at least
two more cusps than does the putative child $N$, then there are
potentially infinitely many Dehn filling
coefficients with normalized length at most
a given bound. Getting around this difficulty
is work in progress.

The overriding difficulty with any Dehn parental
test, regardless of the number of cusps, is that
it necessarily will depend upon
an isometry test for hyperbolic 3-manifolds.
There is an isometry algorithm, due to Scott and
Short, and depending on the hyperbolic structure algorithm of
Manning (see \cite{ScoSho}, \cite{Manning}).
Unfortunately, this algorithm is slow
(though it is elementary-recursive; see \cite{Kup}).
SnapPy therefore does not use such an algorithm,
and instead relies upon a procedure which
is usually fast, but which is not known to
admit a termination argument. In practice,
an implementation of a Dehn parental test
would likely rely on such a heuristic method,
and therefore not be an algorithm. The
existence of an algorithm with running times
similar to those of SnapPy's procedure is
an open question.

If $M$ has exactly one more cusp than $N$, then
one could improve greatly upon the bounds given here by
making rigorous the bounds given in Remark 2.3.b of
\cite{HM}. The number of slopes one must check using
their bounds is $O(\sqrt{1/\Delta V})$, whereas our
bounds are $O(1/\Delta V)$. A rigorous extension
of their bounds to the generic case
$|\pi_0(\partial M)| - |\pi_0(\partial N)| > 1$ could
vastly decrease running time.

The reader may wonder why, after the systole length
test, we don't find the systole, drill it out,
and just test whether or not the result is $M$.
We remind the reader that, as discussed in the
introduction, there is no known, implemented way to
combinatorially represent and drill out given isotopy
classes of geodesics in hyperbolic 3-manifolds.
The point of this paper is that one can avoid drilling
out geodesics. One only needs a nonzero lower bound on
the systole length of $N$.

Developing an algorithm to calculate this systole
length bound is ongoing joint work with Matthias
G\"orner, Neil Hoffman, and Maria Trnkov\'{a}.

\subsection{Applications}\label{sub:applications}
Once the Dehn parental test is finished,
one will be able to calculate
the complexity of 3-manifolds
for certain notions of complexity,
among which is Gabai, Meyerhoff,
and Milley's Mom-number $m$
(see \cite{GMM11}).
The number $m$ has the following nice properties:
\begin{itemize}
\item If $M$ is a Dehn filling of $N$, then $m(M) \leq m(N)$.
\item For $0 < B < \infty$, there is a finite, computable set 
      $S_m(B)$ such that if $m(N) < B$, then $N$ is a Dehn
      filling of some element of $S_m$.
\end{itemize}
The volume function $v$ (by Theorems 3.4 and 3.5 of \cite{Th82})
is known to satisfy these properties as well. However,
whereas $S_m$ is quite easy to implement in software
and relatively quick to execute,
ongoing joint work with T. Crawford,
D. Gabai, R. Meyerhoff, N. Thurston,
and A. Yarmola suggests that no such implementation of an
analog $S_v$ is possible.

Finally, one can rephrase a Dehn parental test as a test
for (directed) adjacency in the big Dehn surgery
graph (see \cite{HofWal15}). An implementation
of a Dehn parental test will enable the construction
of induced subgraphs of the big Dehn surgery graph.

\section*{Acknowledgments}
Martin Bridgeman suggested that
Hodgson and Kerckhoff's work could be
used to develop a Dehn parental test.
Craig Hodgson pointed out that
the notion of normalized length used
in \cite{HK08} has
a precursor in the asymptotics of \cite{NZ85}.
Neil Hoffman pointed out an error in
an earlier draft of this paper.
I had helpful conversations with
Andrew Yarmola and Eric Towers.
The referees made excellent suggestions
for improvements to this paper.
I thank them all.

\bibliography{haraway}{}
\bibliographystyle{amsplain}

\end{document}